\documentclass[11pt]{article}

\usepackage[margin=3cm,a4paper]{geometry}
\usepackage{hyperref}

\hypersetup{
    colorlinks=true,
    linkcolor=blue,
    filecolor=black,
    citecolor=blue,
    urlcolor=blue,
    pdftitle={Extending Exact Integrality Gap Computations for the Metric TSP},
    }

\usepackage{tabularx}
\usepackage{caption}
\usepackage{amssymb,amsmath,amsthm}
\usepackage{mathtools}
\usepackage[capitalise]{cleveref}

\newcommand{\gapp}{\ifmmode\mathrm{gap}^{+}\else{$\mathrm{gap}^+$}\fi}
\newcommand{\nauty}{\texttt{nauty}}

\usepackage[dvipsnames]{xcolor}

\title{Extending Exact Integrality Gap Computations for the Metric TSP}
\author{William Cook\thanks{Combinatorics and Optimization, University of Waterloo, Canada}, Stefan Hougardy\thanks{Research Institute for Discrete Mathematics and Hausdorff Center for Mathematics, University of Bonn, Germany}, Moritz Petrich\thanks{Research Institute for Discrete Mathematics, University of Bonn, Germany}}
\date{\today}

\begin{document}

\maketitle 

\begin{abstract}
The subtour relaxation of the traveling salesman problem (TSP) plays a central role in approximation algorithms and polyhedral studies of the TSP. A long-standing conjecture asserts that the integrality gap of the subtour relaxation for the metric TSP is exactly \(4/3\). In this paper, we extend the exact verification of this conjecture to a larger number of vertices.

Using the framework introduced by Benoit and Boyd in 2008, we confirm their results up to \(n=10\). We further show that for \(n=11\) and \(n=12\), the published lists of extreme points of the subtour polytope are incomplete: one extreme point is missing for \(n=11\) and twenty-two extreme points are missing for \(n=12\). We extend the enumeration of the extreme points of the subtour polytope to instances with up to \(n=16\) vertices in the general case. Restricted to half-integral extreme points, 
we extend the enumeration to $n=18$. Our results provide additional support for the \(4/3\)-Conjecture.

Our \href{https://doi.org/10.60507/FK2/JK95PC}{lists of extreme points} are available on the public bonndata repository.
\end{abstract}

\section{Introduction}
The traveling salesman problem (TSP) is one of the most intensively studied problems in combinatorial 
optimization. Given a finite set of vertices together with pairwise distances, the objective is to determine 
a shortest possible tour that visits each vertex exactly once and returns to the starting point. The problem 
is NP-hard, which motivates the study of polynomial-time approximation algorithms. For the metric 
TSP, the long-standing best-known guarantee was Christofides’  $3/2$-approximation algorithm; this bound was improved to $3/2-\epsilon$ for some fixed $\epsilon>0$ by Karlin, Klein, and Oveis Gharan~\cite{KKO2023}.

A central tool in the analysis of the TSP is the subtour relaxation, which arises from the standard integer-linear-programming formulation of the TSP by imposing the degree constraints and subtour-elimination constraints while dropping integrality requirements. The optimal value of this relaxation provides a natural lower bound on the length of an optimal TSP tour and forms the basis of many approximation algorithms as well as extensive polyhedral investigations.

A long-standing open problem in this context is the $4/3$-Conjecture~\cite{Goe1995,Goe2012,Wil1990}, which asserts that the integrality gap of the subtour relaxation for the metric TSP is at most $4/3$; that is, the ratio between the length of an optimal TSP tour and the optimal value of the subtour relaxation does not exceed $4/3$.
Explicit constructions are known that prove the integrality gap is at least $4/3$~\cite{ Wil1990,Hou2014,HZ2021,Zho2025}, thus showing that the conjectured bound, if true, would be tight. Establishing or refuting the 4/3-Conjecture remains one of the central challenges in the polyhedral study of the TSP.

The $4/3$-Conjecture has been confirmed for several special classes of instances~\cite{BSSS2014,JKW2025,villa2025integralitygaptravelingsalesman}, providing partial evidence for its validity. In addition, it has been verified computationally for small numbers of vertices. In particular, the conjecture has been verified for all instances up to $n=10$ vertices in~\cite{BenoitBoyd2008}, and later extended to all instances up to $n=12$ vertices and to all half-integral extreme points up to $n=14$ in~\cite{BoydElliottMagwood2010,webpageboyd}.

In this work, using the same methods as in~\cite{BenoitBoyd2008}, we further extend these computational results. We verify the conjecture for the general case up to $n=16$ and for the half-integral case up to 
$n=18$. Moreover, we show that the previously published lists of extreme points~\cite{BoydElliottMagwood2010,webpageboyd} are incomplete: 1 extreme point is missing for 
$n=11$ and 22 extreme points are missing for $n=12$. In the half-integral setting, 23 extreme points are missing  in~\cite{BenoitBoyd2008}\footnote{Note that the numbers in Table~2 of~\cite{BenoitBoyd2008} seem not to include the integer tour vertex, contrary to the statement two lines above
the table on page~930.} for $n=13$, $14$ . Our 
\href{https://doi.org/10.60507/FK2/JK95PC}{lists of extreme points} are available on the public bonndata repository~\cite{CHMdata2026}.

\section{Computational Approach}
Benoit and Boyd~\cite{BenoitBoyd2008} described a method to enumerate all extreme points of the subtour polytope for a given number $n$ of vertices. With this approach, they were able to list all extreme points up to $n=10$ within four days of computation time on a SUNW UltraSPARC-II with a single 400MHz processor. They write in their paper~\cite{BenoitBoyd2008}: 
\begin{quote} \textit{
    [$\ldots$] however, for $n=11$, we found that PORTA was unable to generate all of the vertices for $S^n$ using this method, even after running for weeks.}
\end{quote}
In a later paper, Boyd and Elliott-Magwood~\cite{BoydElliottMagwood2010} describe a different 
method that allows them to enumerate all extreme points of the subtour polytope up to $n=12$. With this method, the computation time for $n=11$ was 20 hours on the same computer as before. 

We closely follow the first approach of Benoit and Boyd~\cite{BenoitBoyd2008} but use
stronger theoretical results in Step 2'' described below. Instead of PORTA we use PPL~\cite{Bagnara2008}
to enumerate all extreme points of a polytope. Moreover, we use more powerful hardware: while~\cite{BenoitBoyd2008} used a SUNW UltraSPARC-II with a single 400MHz processor, 
we use a 2.55GHz {AMD EPYC 9684X} machine with 96 cores.
For $n=10$, our hardware/software platform is more than $10^6$ times faster than~\cite{BenoitBoyd2008}.

To enumerate all extreme points of the subtour polytope, we apply the framework described by Benoit and Boyd~\cite{BenoitBoyd2008}. We use the following steps 1''--5'' for each $n$ 
that are very similar to the steps 1'--5' described in~\cite{BenoitBoyd2008}. We start with an empty set 
$\mathcal{Q}^n$ of representatives of non-isomorphic extreme points of the subtour polytope.
\begin{enumerate}
    \item[Step 1''] (\textbf{List potential support graphs}) Generate all non-isomorphic 2-vertex-connected  graphs on $n$ vertices with at most $2n-3$ edges
          and each vertex having degree at least~3. For this step, we use \nauty\ version 2.9.3 \cite{McKayPiperno2014}. 
          This step is identical to Step~1 in~\cite{BoydElliottMagwood2010}.
Results and running times are shown in \cref{tab:support_graphs}. 
\begin{table}[htbp]
    \centering
    \caption{Graphs counted by calling \texttt{callgeng2 -N384 $n$ -Cd3u 0:[$2n-3$]}.}
    \label{tab:support_graphs}
    \renewcommand{\arraystretch}{1.2}
    \begin{tabular}{r@{\hskip 15pt}r@{\hskip 15pt}r}
        $n$ & Number of graphs & Running time in seconds \\
        \hline
         6 & 2 & 0 \\
         7 & 4 & 0 \\
         8 & 38 & 0 \\
         9 & 302 & 0 \\
        10 & 3,745 & 0 \\
        11 & 54,721 & 0 \\
        12 & 956,444 & 0 \\
        13 & 18,957,450 & 0 \\
        14 & 419,857,629 & 8 \\
        15 & 10,244,082,421 & 158 \\
        16 & 272,892,743,415 & 2,371 \\
        \hline
    \end{tabular}
\end{table}

    \item[Step 2''] (\textbf{Enumerate extreme points for support graphs}) 
    For each graph from Step~1'' we check further properties which are necessary for support graphs of extreme points.
    Using a density criterion and 1-block toughness, the number of potential support graphs can be reduced by 25\% for $n=12$ as observed in~\cite{BoydElliottMagwood2010}. 
    New checks are able to remove 96\% of the graphs for $11\leq n\leq16$ safely \cite{Pet2026}.

    For each passed graph $G$ we restrict the subtour polytope by setting all variables corresponding to edges not in $E(G)$ to 0. 
    This yields a face of the subtour polytope or it yields the empty set. 
    For this smaller polytope we enumerate all extreme points with the {Parma Polyhedra Library} (PPL) version 1.2 \cite{Bagnara2008}. 
    Each such generated extreme point, which additionally has $G$ as its support graph\footnote{For example all Hamiltonian cycles of $G$ will also be generated as extreme points, but these are dropped, because their support graph is a strict subgraph of $G$ (unless $G$ itself is a cycle).}, is added to $\mathcal{Q}^n$. 
    Running times of this step dominate those shown in the column ``PPL'' of \cref{tab:running_times}. 

    \item[Step 3''] (\textbf{Generate extreme points via edge splitting operations}) 
    Generate additional extreme points of the subtour polytope by performing edge splitting operations defined in~\cite{BoydElliottMagwood2010} on extreme points in $\mathcal{Q}^{n-1}$. This produces among others all extreme points with adjacent 1-edges. The results and running times for this step are shown in \cref{tab:edgesplitting}.
\begin{table}[ht]
    \centering
    \caption{Non-isomorphic additional extreme points generated exclusively in Step 3''.}
    \label{tab:edgesplitting}
    \renewcommand{\arraystretch}{1.2}
    \begin{tabular}{r@{\hskip 15pt}r@{\hskip 15pt}r}
        $n$ & Number of extreme points & Running time in seconds \\
        \hline
         6 &             1 &      0 \\
         7 &             3 &      0 \\
         8 &             8 &      0 \\
         9 &            40 &      0 \\
        10 &           321 &      0 \\
        11 &         3,616 &      0 \\
        12 &        48,002 &      0 \\
        13 &       738,530 &      1 \\
        14 &    12,188,043 &     27 \\
        15 &   216,746,143 &    690 \\
        16 & 4,007,938,126 & 18,216 \\
        \hline
    \end{tabular}
\end{table}

    \item[Step 4''] (\textbf{Remove isomorphic extreme points})
    Before adding an extreme point to $\mathcal{Q}^n$ we determine a unique representative of its isomorphism class, which may already be contained in $\mathcal{Q}^n$. Otherwise, we add this representative instead of the extreme point. 
    To get the unique representative for an extreme point we apply a canonical labeling of an edge-colored support graph. 
    Each edge is colored with respect to the position of its weight in a sorted list of all the weights that occur in this extreme point.
    Edges with the same weight get identical colors.
    We compute the canonical labeling of the edge colored graph with \nauty\ via a transformation into a layered graph proposed in the  \nauty\  user guide.
    Note that this step is immediately performed inside Step 2'' and Step 3''.

    \item[Step 5''] (\textbf{Compute gap of non-isomorphic extreme points})
    We call an integral extreme point of the subtour polytope a tour.
    For each $x\in \mathcal{Q}^n$ we solve the following LP using IBM ILOG CPLEX Optimization Studio 20.1.0. 
  \begin{align*} 
    \frac{1}{\mathrm{gap}^+(x)} \ \coloneqq\ \min  x^Tc& && \\
    \text{s.t. } c_{\{i,j\}}+c_{\{j,k\}} &\geq c_{\{i,k\}} && \forall i,j,k\in V \\
     c(t) &\geq 1 &&\forall \text{ tours } t \\
    c&\geq 0 &&
  \end{align*}
    The exponentially many tour constraints are handled via a cutting-plane method.
  Floating-point primal and dual solutions are rounded to rational values for which feasibility and optimality are checked exactly. Running times for this step are shown in the column 
  ``gap$^+$'' of \cref{tab:running_times}.
\end{enumerate}
\begin{table}[htbp]
    \centering
    \caption{Running times in seconds for enumerating extreme points (Step~1''--4'', column ``PPL'') and computing the integrality gap (Step~5'', column ``gap$^+$''). The columns ``lrs'' and ``normaliz'' show the runtimes for the general case if instead of PPL we use lrslib or Normaliz.}
    \label{tab:running_times}
    \renewcommand{\arraystretch}{1.2}
    \begin{tabular}{r@{\hskip 15pt}r@{\hskip 15pt}r@{\hskip 15pt}r@{\hskip 15pt}r@{\hskip 15pt}r@{\hskip 15pt}r@{\hskip 15pt}r}
            & \multicolumn{2}{c}{general case} & \multicolumn{2}{c}{half integral case} & & \\
        $n$ & PPL & \gapp & PPL & \gapp & lrs & normaliz \\
        \hline
        11 &         0 &       0 &       0 &      0 &       1 &       0 \\
        12 &         2 &       1 &       0 &      0 &      20 &       2 \\
        13 &        32 &      20 &       0 &      1 &     908 &      34 \\
        14 &       849 &     513 &       6 &      4 & 133,148 &   1,433 \\
        15 &    27,989 &  14,274 &      64 &     29 &       - & 116,936 \\
        16 & 2,107,904 & 432,949 &     656 &    295 &       - &       - \\
        17 &         - &       - &   9,568 &  2,968 &       - &       - \\
        18 &         - &       - & 195,586 & 35,097 &       - &       - \\
        \hline
    \end{tabular}
\end{table}

To enumerate all half-integral extreme points, we modify Step~1'' to generate only graphs with maximum degree~4 and perform Step~2'' similar to~\cite{BenoitBoyd2008}. Running times are shown in \cref{tab:running_times}
below the ``half integral case'' column.

\section{Results}
\cref{tab:vertices_gap} shows the number of extreme points we computed for small values of $n$ and, in addition, the total running times for steps~1''--5''.
Our results match those reported in \cite{BenoitBoyd2008,BoydElliottMagwood2010} up to $n=10$. 
We found 1 and 22, respectively, additional extreme points for $n=11,12$ compared to \cite{webpageboyd}. 
We were able to compute all extreme points up to $n=16$ and all half-integral extreme points up to $n=18$. 
This revealed 5 and 18, respectively, half-integral extreme points more than in \cite{BenoitBoyd2008} for $n=13,14$.
The integrality gap up to $n=16$ is still uniquely attained at the extreme points
 conjectured by Benoit and Boyd~\cite{BenoitBoyd2008}. 
The same holds true when restricted to half-integral extreme points up to $n=18$.

\begin{table}[ht]
    \centering
    \caption{The number of (half-integral) extreme points of the subtour polytope and the maximum integrality gap for a given $n$. Numbers in bold correct or extend the results from~\cite{BenoitBoyd2008, BoydElliottMagwood2010}. Running times are given in seconds.}
    \label{tab:vertices_gap}
    \renewcommand{\arraystretch}{1.2}
    \begin{tabular}{r@{\hskip 20pt}rr@{\hskip 20pt}c@{\hskip 20pt}rr}
        $n$ & \multicolumn{2}{c}{Extreme points} & Integrality gap & \multicolumn{2}{c}{Total running time} \\
        & general & half-integral & (half-integral) & general & half-integral \\
        \hline
         6 &                      2 &                    2 & ${10}/{9} \approx 1.111$    &       0 &       0 \\
         7 &                      3 &                    3 & ${9}/{8} = 1.125$           &       0 &       0 \\
         8 &                     13 &                   12 & ${8}/{7} \approx 1.143$     &       0 &       0 \\
         9 &                     56 &                   42 & ${7}/{6} \approx 1.167$     &       0 &       0 \\
        10 &                    462 &                  208 & ${20}/{17} \approx 1.176$   &       0 &       0 \\
        11 & \textbf{        4,973} &                1,023 & ${19}/{16} \approx 1.188$   &       0 &       0 \\
        12 & \textbf{       68,342} &                5,638 & ${6}/{5} = 1.200$           &       3 &       0 \\
        13 & \textbf{    1,050,837} & \textbf{     31,692} & \boldmath${35}/{29} \approx 1.207$   &        52 &      1 \\
        14 & \textbf{   17,672,908} & \textbf{    185,644} & \boldmath${17}/{14} \approx 1.214$   &     1,362 &     10 \\
        15 & \textbf{  319,185,400} & \textbf{  1,109,906} & \boldmath${11}/{9} \approx 1.222$    &    42,263 &     93 \\
        16 & \textbf{6,165,914,700} & \textbf{  6,780,557} & \boldmath${27}/{22} \approx 1.227$   & 2,540,853 &    951 \\
        17 &                      - & \textbf{ 42,134,667} & \boldmath\textbf{(${53}/{43} \approx 1.233$)} &       - &  12,536 \\
        18 &                      - & \textbf{266,146,869} & \boldmath\textbf{(${26}/{21} \approx 1.238$)} &       - & 230,683 \\
        \hline
    \end{tabular}
\end{table}

All computations were performed on a 2.55GHz {AMD EPYC 9684X} machine with 96 cores.
To verify our results, we checked all computed extreme points with two independent extreme-point tests and one independent isomorphism test.
To provide additional evidence, we also replaced the extreme-point enumeration package PPL by \texttt{lrslib} version 7.3a \cite{Avis2000} and \texttt{Normaliz} version 3.11.1 \cite{Normaliz}.
These packages are based on algorithms other than the double description algorithm used by PPL. 
While \texttt{lrslib} implements a reverse search algorithm, \texttt{Normaliz} is computing with their own pyramid decomposition \cite{BrunsWinfried2016Tpop}.
Due to increased running times reported in the columns ``lrs'' and ``normaliz'' of \cref{tab:running_times} we were able to do this verification only up to $n=14$, respectively $n=15$.

\bibliographystyle{plainurl}
\bibliography{ref}
\end{document}